\documentclass[a4paper,11pt]{article}

\usepackage{bm}
\usepackage{enumerate}
\usepackage{graphicx}
\usepackage{multirow}
\usepackage{multicol}
\usepackage{multirow}
\usepackage{amstext}
\usepackage{amssymb}
\usepackage{pdflscape}
\usepackage{amsfonts}
\usepackage{amsmath}
\usepackage{booktabs}
\usepackage[T1]{fontenc}
\usepackage{lmodern}
\usepackage{natbib}
\usepackage{tikz}
 \bibpunct[, ]{(}{)}{,}{a}{}{,}%
\newtheorem{remark}{Remark}

\renewcommand{\vec}[1]{\boldsymbol{#1}}
\newcommand{\field}[1]{\mathbb{#1}}  % the font for a mathematical field is blackboard
\newcommand{\R}{\field{R}} % the field ofthe uncertainty in  the reals
\newcommand{\Z}{\field{Z}} % the field of the integers
\newcommand{\Uset}{\mathcal{U}} % the perturbation set U
\newcommand{\Zset}{\mathcal{Z}} % the primitive uncertainty set Z
\newcommand{\transp}{^{\top}}
\newcommand{\onenorm}[1]{\left|\left|#1\right|\right|_1}
\newcommand{\twonorm}[1]{\left|\left|#1\right|\right|_2}
\newcommand{\infnorm}[1]{\left|\left|#1\right|\right|_\infty}

\def\Cov{\mathop{\rm Cov}}
\def\Var{\mathop{\rm Var}}

\hypersetup{
    colorlinks=true,
    citecolor=black,
    urlcolor=blue,
    linkcolor=black,
    pdfborder={0 0 0},
}

\begin{document}
	
\title{A Practical Guide to Robust Optimization}

\author{Bram L. Gorissen, {I}hsan Yan{\i}ko\u{g}lu, Dick den Hertog \\ \\
	\textit{\small Tilburg University, Department of Econometrics and Operations Research, 5000 LE Tilburg, Netherlands} \\
	\textit{\small {\tt \{b.l.gorissen,d.denhertog\}@tilburguniversity.edu}, {\tt ihsan.yanikoglu@ozyegin.edu.tr}}}
\date{}
\maketitle

\begin{abstract}
Robust optimization is a young and active research field that has been mainly developed in the last
15 years. Robust optimization is very useful for practice, since it is tailored to the information at hand, and
it leads to computationally tractable formulations. It is therefore remarkable that real-life applications of robust optimization are still lagging behind; there is much more potential
for real-life applications than has been exploited hitherto. The aim of this paper is to help practitioners to understand robust optimization and to
successfully apply it in practice. We provide a brief introduction to robust optimization, and also describe important do's and don'ts for using it in practice.
We use many small examples to illustrate our discussions.
\end{abstract}

\begin{tikzpicture}[remember picture,overlay]
\node[anchor=south,yshift=10pt] at (current page.south) {\fbox{\parbox{\dimexpr\textwidth-\fboxsep-\fboxrule\relax}{\footnotesize This is an author-created, un-copyedited version of an article published in Omega \href{http://dx.doi.org/10.1016/j.omega.2014.12.006}{DOI:10.1016/j.omega.2014.12.006}.}}};
\end{tikzpicture}

\section{Introduction}\label{II}
Real-life optimization problems often contain uncertain data. Data can be inherently stochastic/random or it can be uncertain due to errors. The reasons for data errors could be measurement/estimation errors that come from the lack of knowledge of the parameters of the mathematical model (e.g., the uncertain demand in an inventory model) or could be implementation errors that come from the physical impossibility to exactly implement a computed solution in a real-life setting. There are two approaches to deal with data uncertainty in optimization, namely robust and stochastic optimization. \emph{Stochastic optimization} (SO) has an important assumption, i.e., the true probability distribution of uncertain data has to be known or estimated. If this condition is met and the reformulation of the uncertain optimization problem is computationally tractable, then SO is the methodology to solve the uncertain optimization problem at hand. For details on SO, we refer to \citet{pre,birge} and \citet{shap}, but the list of references can be easily extended.

\emph{Robust optimization} (RO), on the other hand, does not assume that probability distributions are known, but instead it assumes that the uncertain data resides in a so-called uncertainty set. Additionally, basic versions of RO assume ``hard'' constraints, i.e., constraint violation cannot be allowed for any realization of the data in the uncertainty set. RO is popular because of its computational tractability for many classes of uncertainty sets and problem types. For a detailed overview of the RO framework, we refer to \citet{BenTal,benem} and \citet{bert}.

Although the first published study \citep{soyster} dates back to 1970s, RO is a relatively young and active research field, and has been mainly developed in the last 15 years. There have been many publications that show the value of RO in many fields of application including finance \citep{finance}, energy \citep{bertsimas2013adaptive, energy}, supply chain \citep{BenTal:2005:retailer,lim2013joint}, healthcare \citep{health}, engineering \citep{engineer}, scheduling \citep{yan2009inter}, marketing \citep{wang2012robust}, etc. Indeed, the RO concepts and techniques are very useful for practice, since they are tailored to the information at hand and leads to tractable formulations. It is therefore remarkable that real-life applications are still lagging behind; there is much more potential for real-life applications than has been exploited hitherto.

In this paper we give a concise description of the basics of RO, including so-called {\it adjustable} RO for multi-stage optimization problems.
Moreover, we extensively discuss  several items that are important when applying RO, and that are often not well understood or incorrectly applied by practitioners.
Several important do's and don'ts are discussed, which may help the practitioner to successfully apply RO.
We use many small examples to illustrate our discussions.

The remainder of the paper is organized as follows.
Section \ref{DR} gives a concise introduction to RO.
Section \ref{B1}-\ref{B4} discuss several important practical issues in more detail.
An important ingredient of RO is the so-called {\it uncertainty set}, which is the set of values for the uncertain parameters that are taken into account in the robust problem. This set has to be specified by the user, and Section \ref{B1} presents several ways how to construct uncertainty sets.
Section \ref{B3} discusses an important technical aspect in {\it adjustable} RO.
In practice, multi-stage problems may contain adjustable variables that are integer. Section \ref{I1} proposes a RO method that deals with such integer variables.
Section \ref{D2} warns that the robust versions of equivalent deterministic problems are not necessarily equivalent, and gives the correct robust counterparts. Section \ref{D1} discusses that equality constraints containing uncertain parameters needs a special treatment. In RO one optimizes for the worst case scenario, i.e., for the worst case values of the uncertain parameters. Although this statement is true, the way it is formulated is often misunderstood. In Section \ref{I2}, we therefore clarify how it should be formulated. It is important to test how robust the final solution is, and to compare it to, e.g., the nominal solution. Section \ref{B2} discusses how to assess the robustness performance of a given solution via a simulation study. Section \ref{B4} shows that RO applied in a folding horizon setting may yield better solutions than {\it adjustable} RO for multi-stage problems.
Section \ref{conc} summarizes our conclusions, and indicates future research topics.

\section{Introduction to robust optimization}\label{DR}
In this section we first give a brief introduction to RO, and then we give a procedure for applying RO in practice. The scopes of the following sections are also presented in this section.
\subsection{Robust optimization paradigm}
For the sake of exposition, we use an uncertain linear optimization problem, but we point out that most of our discussions in this paper can be generalized for other classes of uncertain optimization problems. The ``general'' formulation of the uncertain linear optimization problem is as follows:
\begin{equation}
\min_{\vec{x}}\; \{\vec{c \transp x} : \vec{Ax} \leq \vec{d}\}_{(\vec{c},\vec{A},\vec{d}) \in \Uset }, \label{eq:generalrc}
\end{equation}
where $\vec{c} \in \R^n$, $\vec{A}\in \R^{m \times n}$ and $\vec{d}\in \R^m$ denote the uncertain coefficients, and $\Uset$ denotes the user specified uncertainty set. The ``basic'' RO paradigm is based on the following three assumptions \citep[p.~xii]{BenTal}:
\medskip
\begin{enumerate}
\item[\textbf{A.1.}] All decision variables $\vec{x}\in \R^n$ represent ``here and now'' decisions: they should get specific numerical values as a result of solving the problem before the actual data ``reveals itself''.
\item[\textbf{A.2.}] The decision maker is fully responsible for consequences of the decisions to be made when, and only when, the actual data is within the prespecified uncertainty set $\Uset$.
\item[\textbf{A.3.}] The constraints of the uncertain problem in question are ``hard'', i.e., the decision maker cannot tolerate violations of constraints when the data is in the prespecified uncertainty set $\Uset$.
\end{enumerate}
\medskip
In addition to the ``basic'' assumptions, we may assume without loss of generality that: 1) the objective is certain; 2) the constraint right-hand side is certain; 3) $\Uset$ is compact and convex; and 4) the uncertainty is constraint-wise. Below, we explain the technical reasons of why these four assumptions are not restrictive.
\medskip
\begin{enumerate}
  \item[\textbf{E.1.}] Suppose the objective coefficients ($\vec{c}$) are uncertain and (say) these coefficients reside in the uncertainty set $\mathcal{C}$:
  \[\min_{\vec{x}} \ \max_{\vec{c}\in \mathcal{C}} \ \{\vec{c}^\top \vec{x} : \vec{A}\vec{x} \leq \vec{d} \quad \forall \vec{A} \in \Uset\}.\]
  Without loss of generality we may assume that the uncertain objective of the optimization problem can be equivalently reformulated as certain \citep[p.~10]{BenTal}:
  \[\min_{\vec{x},\ t} \ \{t : \vec{c}^\top \vec{x} - t \leq 0 \quad \forall \vec{c} \in \mathcal{C}, \ \vec{A}\vec{x} \leq \vec{d} \quad \forall \vec{A} \in \Uset\}, \]
  using a reformulation and the additional variable $t\in \mathbb{R}$.
  \item[\textbf{E.2.}] The second assumption is not restrictive because the uncertain right-hand side of a constraint can always be translated to an uncertain coefficient by introducing an extra variable $x_{n+1}=-1$.
  \item[\textbf{E.3.}] The uncertainty set $\Uset$ can be replaced by its convex hull conv($\Uset$), i.e., the smallest convex set that includes $\Uset$, because testing the feasibility of a solution with respect to $\Uset$ is equivalent to taking the supremum of the left hand side of a constraint over $\Uset$, which yields the same optimal objective value if the maximization is conv($\Uset$). For details of the formal proof and the compactness assumption, see \citep[pp.~12--13]{BenTal}.
  \item[\textbf{E.4.}] To illustrate that robustness with respect to $\Uset$ can always be formulated constraint-wise, consider a problem with two constraints and with uncertain parameters $d_1$ and $d_2$: $x_1 + d_1 \leq 0$, $x_2 + d_2 \leq 0$. Let $\Uset=\{\vec{d} \in\mathbb{R}^2 : d_1\geq 0, \ d_2 \geq 0, \ d_1 + d_2 \leq 1\}$ be the uncertainty set. Then, $\Uset_i = [0,1]$ is the projection of $\Uset$ on $d_i$. It is easy to see that robustness of the $i$-th constraint with respect to $\Uset$ is equivalent to robustness with respect to $\Uset_i$, i.e., the uncertainty in the problem data can be modelled \emph{constraint-wise}. For the general proof, see \citep[pp.~11--12]{BenTal}.
\end{enumerate}
\medskip
For uncertain \emph{nonlinear} optimization problems, the assumptions are also without loss of generality, except the third basic assumption [E.3].

If we assume $\vec{c}\in \R^n$ and $\vec{d}\in \R^m$ are certain, then the robust reformulation of \eqref{eq:generalrc} that is generally referred to as the \emph{robust counterpart} (RC) problem is given as follows:
\begin{equation}
\min_{\vec{x}} \{\vec{c \transp x}: \;\; \vec{A(\bm{\zeta})x} \leq \vec{d} \quad \forall \bm{\zeta} \in \Zset \big\}, \label{eq:si}
\end{equation}
where $\Zset \subset \R^L$ denotes the user specified \emph{primitive} uncertainty set. A solution $\vec{x}\in \R^n$ is called \emph{robust feasible} if it satisfies the uncertain constraints [$\vec{A(\bm{\zeta})x} \leq \vec{d}$] for all realizations of $\bm{\zeta}\in \Zset$.

As it is mentioned above, and explained in [E.4], we may focus on a single constraint, since the uncertainty is constraint-wise in RO. A single constraint taken out of \eqref{eq:si} can be modeled as follows:
\begin{equation}\label{s0}
\vec{(a+P\bm{\zeta}) \transp x} \leq d  \quad \forall \bm{\zeta} \in \Zset.
\end{equation}
In the left-hand side of \eqref{s0}, we use a factor model to formulate a single constraint of \eqref{eq:si} as an affine function $\vec{a+P}\bm{\zeta}$ of the \emph{primitive} uncertain parameter $\bm{\zeta} \in \Zset$, where $\vec{a} \in \R^n$ and $\vec{P}\in \R^{n \times L}$. One of the most famous example of such a factor model is the 3-factor model of \citet{fama}, which models different type of assets as linear functions of a limited number of uncertain economic factors. To point out, the dimension of the general uncertain parameter $\vec{P}\bm{\zeta}$ is often much higher than that of the primitive uncertain parameter $\bm{\zeta}$ (i.e., $n \gg L$).

\subsection{Solving the robust counterpart}\label{sec:solvingrc}
Notice that \eqref{s0} contains infinitely many constraints due to the \emph{for all} ($\forall$) quantifier imposed by the worst case formulation, i.e., it seems intractable in its current form. There are two ways to deal with this. The first way is to apply robust reformulation techniques to exclude the \emph{for all} ($\forall$) quantifier. If deriving such a robust reformulation is not possible, then the second way is to apply the adversarial approach. In this section, we describe the details of these two approaches.

We start with the first approach, which consists of three steps. The result will be a computationally tractable RC of \eqref{s0}, which contains a finite number of constraints. Note that this reformulation technique is one of the main techniques in RO \citep{bert}.

We illustrate the three steps of deriving the RC based on a polyhedral uncertainty set:
\[\mathcal{Z}=\{\bm{\zeta}:\vec{D\bm{\zeta} + q} \geq \bm{0}\},\]
where $\vec{D}\in \R^{m \times L}, \bm{\zeta}\in \R^L,$ and $\vec{q}\in \R^m$.

\noindent\textbf{Step 1} (\emph{Worst case reformulation}). Notice that \eqref{s0} is equivalent to the following worst case reformulation:
\begin{equation}\label{s1}
\vec{a^\top x} +  \max_{\bm{\zeta}  : \ \vec{D\bm{\zeta} + q} \geq 0} (\vec{P^\top x})^\top \bm{\zeta} \leq d.
\end{equation}

\noindent\textbf{Step 2} (\emph{Duality}). We take the dual of the inner maximization problem in \eqref{s1}. The inner maximization problem and its dual yield the same optimal objective value by strong duality. Therefore, \eqref{s1} is equivalent to
\begin{equation}\label{s2}
\quad \vec{a^\top x} +  \min_{\vec{w}} \{\vec{q^\top w   : \ D^\top w = -P^\top x, \ w} \geq \vec{0}\} \leq d.
\end{equation}
\textbf{Step 3} (\emph{RC}). It is important to point out that we can omit the minimization term in \eqref{s2}, since it is sufficient that the constraint holds for at least one $\vec{w}$. Hence, the final formulation of the RC becomes
\begin{equation}\label{s3}
\exists \ \vec{ w}: \ \ \vec{a^\top x +  q^\top w} \leq d,  \quad     \vec{D^\top w = -P^\top x} , \quad     \vec{w} \geq \vec{0}.
\end{equation}
Note that the constraints in \eqref{s3} are linear in $\vec{x}\in \R^n$ and $\vec{w}\in \R^m$.

\begin{table}[h]
\centering
\caption{Tractable reformulations for the uncertain constraint [$ \vec{(a+P\bm{\zeta}) \transp x} \leq d  \quad \forall \bm{\zeta} \in \Zset$]  \newline \noindent for different types of uncertainty sets}\label{T1}
\begin{tabular}{l l l l }
\toprule
 Uncertainty &$\Zset$ & Robust Counterpart & Tractability \\
\midrule
Box  & $\|\bm{\zeta}\|_\infty \leq 1$ & $ \vec{a\transp x} + \| \vec{P \transp x} \|_1 \leq d $ & LP\\ %\midrule
Ellipsoidal & $\|\bm{\zeta}\|_2 \leq 1$ & $ \vec{a\transp x} + \|\vec{P \transp x} \|_2 \leq d $ & CQP\\ %\midrule
Polyhedral & $ \vec{D\bm{\zeta}+q} \geq \bm{0}$ & $ \begin{cases} \vec{a \transp x} + \vec{q \transp w} \leq d \\ \vec{D \transp w = -P \transp x} \\ \vec{w}\geq \vec{0}  \end{cases} $  & LP\\ %\midrule
Cone {\tiny (closed, convex, pointed)} & $ \vec{D\bm{\zeta}+q} \in K$ & $ \begin{cases}
    \vec{a^\top x + q^\top w} \leq d \\
   \vec{D^\top w = -P^\top x} \\
   \vec{w}\in K^*
  \end{cases} $& Conic Opt.\\
Convex cons. & $ h_k(\bm{\zeta}) \leq 0 \quad \forall k$ & $ \begin{cases} \vec{a \transp x} + \sum_k u_k h_k^*\left(  \frac{\vec{w^k}}{u_k}\right)  \leq d \\ \sum_k\vec{ w^k = P \transp x} \\ \vec{u} \geq \vec{0} \end{cases} $  &Convex Opt. \\
\bottomrule
\multicolumn{4}{l}{\footnotesize{$(\ast)$ $h^*$ is the convex conjugate function, i.e,  $h^*(x) = \sup_y \{ x^\top y - h(y) \}$; and $K^*$ is the dual cone of $K$;}}\\
\multicolumn{4}{l}{\footnotesize{Table taken from \citet{Bental2011}}}
\end{tabular}
\end{table}

Table \ref{T1} presents the tractable robust counterparts of an uncertain linear optimization problem for different classes of uncertainty sets. These robust counterparts are derived using the three steps that are described above. However, we need conic duality instead of LP duality in Step 2 to derive the tractable robust counterparts for the conic uncertainty set; see the fourth row of Table \ref{T1}. Similarly, to derive the tractable RC for an uncertainty region specified by general convex constraints, i.e., in the fifth row of Table \ref{T1}, we need Fenchel duality in Step 2; see \citet{rocka} for details on Fenchel duality, and \citet{Bental2011} for the formal proof of the associated RC reformulation. Notice that each RC constraint has a positive safeguard in the constraint left-hand side, e.g., $\| \vec{P^\top x} \|_1$, $ \| \vec{P^\top x }\|_2$, and $\vec{q^\top w}$; see the tractable RCs in the third column of Table \ref{T1}. These safeguards represent the level of robustness that we introduce to the constraints. Note the equivalence between the robust counterparts of the polyhedral and the conic uncertainty set when $K$ is the nonnegative orthant.

\emph{Nonlinear problems.} Table \ref{T1} focuses on uncertain linear optimization problems, i.e., both linear in the decision variables and the uncertain parameters. Notice that, different than the presented results, the original uncertain optimization problem can be nonlinear in the optimization variables and/or the uncertain parameters; for more detailed treatment of such problems that are nonlinear in the uncertain parameters, we refer to \citet[pp. 383--388]{BenTal} and \citet{Bental2011}.

\emph{Adversarial approach.} If the robust counterpart cannot be written as or approximated by a tractable reformulation, we advocate to perform the so-called {\it adversarial approach}. The adversarial approach starts with a finite set of scenarios $S_i \subset \Zset_i$ for the uncertain parameter in constraint $i$. For example, at the start, $S_i$ only contains the nominal scenario. Then, the robust optimization problem, which has a finite number of constraints since $\Zset_i$ has been replaced by $S_i$, is solved. If the resulting solution is robust feasible, we have found the robust optimal solution. If that is not the case, we can find a scenario for the uncertain parameter that makes the last found solution infeasible, e.g., we can search for the scenario that maximizes the infeasibility. We add this scenario to $S_i$, and solve the resulting robust optimization problem, and so on. For a more detailed description, we refer to \citet{bienstock}. It appeared that this simple approach often converges to optimality in a few number of iterations. The advantage of this approach is that solving the robust optimization problem with $S_i$ instead of $\Zset_i$ in each iteration, preserves the structure of the original optimization problem. Only constraints of the same type are added, since constraint $i$ should hold for all  scenarios in $S_i$. This approach could be faster than reformulating, even for polyhedral uncertainty sets. See \citet{Bert2014} for a comparison. Alternatively, if the probability distribution of the uncertain parameter is known, one may also use the \emph{randomized sampling} of the uncertainty set proposed by \citet{calafiore}. The randomized approach substitutes the infinitely many robust constraints with a finite set of constraints that are randomly sampled. It is shown that such a randomized approach is an accurate approximation of the original uncertain problem provided that a sufficient number of samples is drawn; see \citet[Theorem 1]{campi}.

\emph{Pareto efficiency.} \citet{Dan} discovered that ``the inherent focus of RO on optimizing performance only under worst case outcomes might leave decisions un-optimized in case a non worst case scenario materialized''. Therefore, the ``classical'' RO framework might lead to a Pareto inefficient solution; i.e., an alternative robust optimal solution may guarantee an improvement in the objective or slack size for (at least) one scenario without deteriorating it in other scenarios. Given a robust optimal solution, Iancu and Trichakis propose optimizing a new problem to find a solution that is Pareto efficient. In this new problem, the objective is optimized for a scenario in the interior of the uncertainty set, e.g., for the nominal scenario, while the worst case objective is constrained to be not worse than the robust optimal objective value. For more details on Pareto efficiency in robust linear optimization we refer to \citet{Dan}.

\subsection{Adjustable robust optimization}\label{sec:adjintro}
In multistage optimization, the first assumption [A.1] of the RO paradigm, i.e., the decisions are ``here and now'', can be relaxed. For example, the amount a factory will produce next month is not a ``here and now'' decision, but a ``wait and see'' decision that will be taken based on the amount sold in the current month. Some decision variables can therefore be adjusted at a later moment in time according to a decision rule, which is a function of (some or all part of) the uncertain data. The \emph{adjustable} RC (ARC) is given as follows:
\begin{align}
\min_{\vec{x, y(\cdot)}} \{\vec{c\transp x}: \;\; \vec{A}(\bm{\zeta})\vec{x} + \vec{B} \vec{y}(\bm{\zeta})\leq \vec{d} \quad \forall \bm{\zeta} \in \Zset \}, \label{adjust}
\end{align}
where $\vec{x}\in \R^n$ is the first-stage ``here and now'' decision that is made before $\bm{\zeta}\in \R^L$ is realized, $\vec{y}\in \R^k$ denotes the second-stage ``wait and see'' decision that can be adjusted according to the actual data, and $\vec{B}\in \R^{m \times k}$ denotes a certain coefficient matrix (i.e., fixed recourse).

However, ARC is a complex problem unless we restrict the function $\vec{y}(\bm{\zeta})$ to specific classes; see \citet[Ch. 14]{BenTal} for details. In practice, $\vec{y}(\bm{\zeta})$ is often approximated by affine (or linear) decision rules:
\begin{equation}
\vec{y}(\bm{\zeta}) := \vec{y^0+Q} \bm{\zeta}, \label{decision}
\end{equation}
because they yield computationally tractable \emph{affinely} ARC (AARC) reformulations, where $\vec{y^0}\in \R^k$ and $\vec{Q} \in \R^{k\times L}$ are the coefficients in the decision rule, which are to be optimized. Eventually, the tractable reformulation of the constraints in \eqref{adjust}:
\[\min_{\vec{x, y^0, Q}}\{ \vec{c^\top x}: \;\; \vec{A(\bm{\zeta})x+By^0+BQ\bm{\zeta}} \leq \bm{d}  \quad \forall \bm{\zeta} \in \Zset \} \]
can be derived in a similar vein by applying the three steps that are described above, since the problem is affine in the uncertain parameter $\bm{\zeta}$, and the decision variables $\vec{x, y^0}$, and $\vec{Q}$.

To point out, ARO is less conservative than the classic RO approach, since it yields more flexible decisions that can be adjusted according to the realized portion of data at a given stage. More precisely, ARO yields optimal objective values that are at least as good as that of the standard RO approach. In addition, aside from introducing additional variables and constraints, the AARC does not introduce additional computational complexity to that of RO with fixed recourse, and it can be straightforwardly adopted to the classic RO framework, i.e., the AARC is a computationally tractable approach. Moreover, affine decision rules are optimal or near optimal for many practical cases; e.g., inventory management \citep{inventory}. Last but not least, ARO has many applications in real-life, e.g., supply chain management \citep{BenTal:2005:retailer}, project management \citep[Ex. 14.2.1]{BenTal}, engineering \citep{engineer}, and so on.

Notice that we adopt affine decision rules in the ARC, but it is important to point out that tractable ARC reformulations for nonlinear decision rules also exist for specific classes; we refer to \citet[Ch.~14.3]{BenTal} and \citet{lift}.

\emph{Integer adjustable variables.} A parametric decision rule, like  the linear one in \eqref{decision}, cannot be used for integer adjustable variables, since we have then to enforce that the decision rule is integer for all $\bm{\zeta} \in \Zset$. In Section \ref{I1} we propose a general way for dealing with adjustable integer variables similar to \citet{Bert2010}. However, much more research is needed.
\vspace{-0.5cm}
\begin{table}
\setlength{\tabcolsep}{9pt}
\caption{}\label{procedure} \begin{tabular}{l p{10.5cm} l} \hline
\multicolumn{2}{l}{\textbf{Practical RO procedure}} &\textbf{Ref. section(s)}\\ \hline \\
\textbf{Step 0:} &Solve the nominal problem. &- \\ \\
\textbf{Step 1:} &\textbf{a)} Determine the uncertain parameters. & \multirow{2}{*}{\S \ 3} \\[4pt]
                 &\textbf{b)} Determine the uncertainty set. &\\ \\
\textbf{Step 2:} &Check robustness of the nominal solution. & \multirow{2}{*}{\S \ 9} \\[4pt]
                 &\textbf{IF} the nominal solution is robust ``enough'' \textbf{THEN} stop. &\\ \\
\textbf{Step 3:} &\textbf{a)} Determine the adjustable variables. & \S\S \ 5 \\[4pt]
                 &\textbf{b)} Determine the type of decision rules. & \S\S \ 4, 10 \\ \\
\textbf{Step 4:} &Formulate the robust counterpart. & \S\S \ 6, 7, 8 \\ \\
\textbf{Step 5:} &Solve the (adjustable) robust counterpart via an exact or approximate tractable reformulation, or via the adversarial approach. &\multirow{2}{*}{\S \ 2}  \\ \\[8pt]
\textbf{Step 6:} &Check quality of the robust solution. &\multirow{2}{*}{\S \ 9} \\[4pt]
                 &\textbf{IF} the solution is ``too conservative'' \textbf{THEN} go to \textbf{Step 1b} or \textbf{Step 3}.  \\ \\ \hline
\end{tabular}
\end{table}

\subsection{Robust optimization procedure}

Now having introduced the general notation in RO and \emph{adjustable} RO (ARO), we can give a procedure for applying RO in practice; see Table \ref{procedure}.

\noindent In the remainder of this paper, we describe the most important items at each step of this procedure. The associated section(s) for each step are reported in the last column of Table \ref{procedure}.

\section{Choosing the uncertainty set}\label{B1}
In this section we describe different possible uncertainty sets and their advantages and disadvantages. Often one wants to make a trade-off between robustness against each physical realization of the uncertain parameter and the size of the uncertainty set. A box uncertainty set that contains the full range of realizations for each component of $\vec{\zeta}$ is the most robust choice and guarantees that the constraint is never violated, but on the other hand there is only a small chance that all uncertain parameters take their worst case values. This has led to the development of smaller uncertainty sets that still guarantee that the constraint is ``almost never'' violated. Such guarantees are inspired by chance constraints, which are constraints that have to hold with at least a certain probability. Often the underlying probability distribution is not known, and one seeks a distributionally robust solution. One application of RO is to provide a tractable safe approximation of the chance constraint in such cases, i.e., a tractable formulation that guarantees that the chance constraint holds:
\begin{align} \mbox{if $\vec{x}$ satisfies } \vec{a}(\vec{\zeta}) \transp \vec{x} \leq d \quad \forall \vec{\zeta} \in \Uset_\varepsilon \mbox{, then $\vec{x}$ also satisfies } \mathbb{P}_{\vec{\zeta}}( \vec{a}(\vec{\zeta}) \transp \vec{x} \leq d ) \geq 1-\varepsilon. \label{def:chanceconstraint} \end{align}
For $\varepsilon = 0$, a chance constraint is a traditional robust constraint. The challenge is to determine the set $\Uset_\varepsilon$ for other values of $\varepsilon$. We distinguish between uncertainty sets for uncertain parameters and for uncertain probability vectors.

For uncertain parameters, many results are given in \citep[Chapter 2]{BenTal}. The simplest case is when the only knowledge about $\vec{\zeta}$ is that $\infnorm{\vec{\zeta}} \leq 1$.  For this case, the box uncertainty set is the only set that can provide a probability guarantee (of $\varepsilon = 0$). When more information becomes available, such as bounds on the mean or variance, or knowledge that the probability distribution is symmetric or unimodal, smaller uncertainty sets become available. \citet[Table 2.3]{BenTal} list seven of these cases. Probability guarantees are only given when $\infnorm{\vec{\zeta}} \leq 1$, $\mathbb{E}(\vec{\zeta}) = \vec{0}$ and the components of $\vec{\zeta}$ are independent. We mention the uncertainty sets that are used in practice when box uncertainty is found to be too pessimistic. The first is an ellipsoid \citep[Proposition 2.3.1]{BenTal}, possibly intersected with a box \citep[Proposition 2.3.3]{BenTal}:
\begin{align}
\Uset_\varepsilon = \{ \vec{\zeta} : \twonorm{\vec{\zeta}} \leq \Omega \quad \infnorm{\vec{\zeta}} \leq 1 \}, \label{unc:bental}
\end{align}
where $\varepsilon = \exp(-\Omega^2 / 2)$. The second is a polyhedral set \citep[Proposition 2.3.4]{BenTal}, called budgeted uncertainty set or the ``Bertsimas and Sim'' uncertainty set \citep{BertsimasPoR}:
\begin{align}
\Uset_\varepsilon = \{ \vec{\zeta} : \onenorm{\vec{\zeta}} \leq \Gamma \quad \infnorm{\vec{\zeta}} \leq 1 \}, \label{unc:bands}
\end{align}
where $\bm{\zeta} \in \mathbb{R}^L$, and $\varepsilon = \exp(-\Gamma^2 / (2L))$. A stronger bound is provided in \citep{BertsimasPoR}. This set has the interpretation that (integer) $\Gamma$ controls the number of elements of $\vec{\zeta}$ that may deviate from their nominal values. \eqref{unc:bental} leads to better objective values for a fixed $\varepsilon$ compared to \eqref{unc:bands}, but gives rise to a CQP for an uncertain LP while \eqref{unc:bands} results in an LP and is therefore more tractable from a computational point of view.

\citet{Bandi2012} propose uncertainty sets based on the central limit theorem. When the components of $\bm{\zeta} \in \mathbb{R}^L$ are independent and identically distributed with mean $\mu$ and variance $\sigma^2$, the uncertainty set is given by:
\begin{align*}
\Uset_\varepsilon = \left\{ \vec{\zeta} : \left| \sum_{i=1}^L \zeta_i - L \mu \right| \leq \rho \sqrt{n} \sigma \right\},
\end{align*}
where $\rho$ controls the probability of constraint violation $1-\varepsilon$. Bandi and Bertsimas also show variations on $\Uset_\varepsilon$ that incorporate correlations, heavy tails, or other distributional information. The advantage of this uncertainty set is its tractability, since the robust counterpart of an LP with this uncertainty set is also LP. A disadvantage of this uncertainty set is that it is unbounded for $L>1$, since one component of $\vec{\zeta}$ can be increased to an arbitrarily large number (while simultaneously decreasing a different component). This may lead to intractability of the robust counterpart or to trivial solutions. In order to avoid infeasibility, it is necessary to define separate uncertainty sets for each constraint, where the summation runs only over the elements of $\vec{\zeta}$ that appear in that constraint. Alternatively, it may help to take the intersection of $\Uset_\varepsilon$ with a box.

\citet{bertsimas2013data} show how to construct uncertainty sets based on historical data and statistical tests. An advantage compared to the aforementioned approach, is that it requires fewer assumptions. For example, independence of the components of $\bm{\zeta}$ is not required. They have included a section with practical recommendations, which is too extensive to discuss here.

\citet{bertsimas2009constructing} show how to construct uncertainty sets based on coherent risk measures that model the preferences of the decision maker, and on past observations of the uncertain data. Starting from a representation theorem of coherent risk measures, they make a natural link to RO. For a specific class of risk measures, this gives rise to polyhedral uncertainty sets.

We now focus on uncertain probability vectors, i.e,. $\Uset \subset \Delta^{L-1} = \{ \vec{p} \in \R^L : \vec{p} \geq \vec{0}, \; \sum_{i=1}^L p_i = 1 \}$. These appear, e.g., in a constraint on expected value or variance. \citet{BTHWMR2011} construct uncertainty sets based on $\phi$-divergence. The $\phi$-divergence between the vectors $\vec{p}$ and $\vec{q}$ is:
\[ I_\phi(\vec{p},\vec{q}) = \sum_{i=1}^L q_i \phi\left( \frac{p_i}{q_i} \right), \]
where $\phi$ is the (convex) $\phi$-divergence function. Let $\vec{p}$ denote a probability vector and let $\vec{q}$ be the vector with observed frequencies when  $N$ items are sampled according to $\vec{p}$. Under certain regularity conditions,
\[ \frac{2N}{\phi^{''}(1)} I_\phi(\vec{p},\vec{q}) \stackrel{d}{\rightarrow} \chi^2_{L-1} \textrm{ as } N \to \infty. \]
This motivates the use of the following uncertainty set:
\[ \Uset_\varepsilon = \{ \vec{p} : \vec{p} \geq \vec{0}, \quad \vec{e} \transp \vec{p} = 1, \quad \frac{2N}{\phi^{''}(1)} I_\phi(\vec{p},\hat{\vec{p}}) \leq \chi^2_{L-1; 1-\varepsilon} \}, \]
where $\hat{\vec{p}}$ is an estimate of $\vec{p}$ based on $N$ observations, and $\chi^2_{L-1; 1-\varepsilon}$ is the $1-\varepsilon$ percentile of the $\chi^2$ distribution with $L-1$ degrees of freedom. The uncertainty set contains the true $\vec{p}$ with (approximate) probability $1-\varepsilon$. \citet{BTHWMR2011} give many examples of $\phi$-divergence functions that lead to tractable robust counterparts.

An alternative to $\phi$-divergence is using the Anderson-Darling test to construct the uncertainty set \citep[Ex. 15]{Bental2011}.

For multistage problems, \citet{goh2010distributionally} show how to construct uncertainty sets based on, e.g., the covariance matrix or bounds on the expected value. They show how these sets can be used for optimizing piecewise-linear decision rules.

We conclude this section by pointing out a mistake that is sometimes made regarding the interpretation of the uncertainty set. Sometimes the set $\Uset$ in \eqref{def:chanceconstraint} is constructed such that it contains the true parameter with probability $\varepsilon$. This provides a much stronger probability guarantee than one expects. For example, let the uncertain parameter $\vec{\zeta} \in \R^{L}$ with $L=10$ have independent components with mean 0 and variance 1. Then, $\vec{\zeta} \transp \vec{\zeta} \sim \chi^2_L$, so the set $\Uset_\varepsilon = \{\vec{\zeta} : \twonorm{\vec{\zeta}} \leq \sqrt{ \chi^2_{L; 1-\varepsilon} } \}$ contains $\vec{\zeta}$ with probability $1-\varepsilon$. Consequently, the logical implication \eqref{def:chanceconstraint} holds. However, the probability in \eqref{def:chanceconstraint} is much larger than $\varepsilon$, since the constraint also holds for the ``good'' realizations of the uncertain parameter outside the uncertainty set. For example, the singleton $\Uset_\varepsilon = \{ \vec{0} \}$ satisfies $\mathbb{P}(\vec{\zeta} \in \Uset_\varepsilon) = 0$, but the probability on the right hand side of \eqref{def:chanceconstraint} is 0.5 if the components of $\vec{\zeta}$ are independent and have a median of $\vec{0}$. To see this, note that $\vec{a}(\vec{\zeta}) \transp \vec{x}$ is an affine function of $\vec{\zeta}$, and therefore has a median of $d = \vec{a}(\vec{0}) \transp \vec{x}$. So, the probability that $\vec{a}(\vec{\zeta}) \transp \vec{x}$ becomes larger than $d$ equals the probability that it becomes smaller. In order to construct the correct set $\Uset_{\varepsilon}$, we first write the explicit chance constraint. Since $(\vec{a} + \vec{P} \vec{\zeta}) \transp \vec{x} \leq d$ is equivalent to $\vec{a} \transp \vec{x} + (\vec{P} \transp \vec{x}) \transp \vec{\zeta} \leq d$, and since the term $(\vec{P} \transp \vec{x}) \transp \vec{\zeta}$ follows a normal distribution with mean 0 and standard deviation $\twonorm{\vec{P} \transp \vec{x}}$, the chance constraint can explicitly be formulated as $\vec{a} \transp \vec{x} + z_{1-\varepsilon} \twonorm{\vec{P} \transp \vec{x}} \leq d$, where $z_{1-\varepsilon}$ is the $1-\varepsilon$ percentile of the normal distribution. This is the robust counterpart of the original linear constraint with ellipsoidal uncertainty and a radius of $z_{1-\varepsilon}$. The value $z_{1-\varepsilon}=9.3$ coincides with $\varepsilon \approx 7.0 \cdot 10^{-21}$. So, while one thinks to construct a set that makes the constraint hold in 50\% of the cases, the set actually makes the constraint hold in almost all cases. To make the chance constraint hold with probability $1-\varepsilon$, the radius of the ellipsoidal uncertainty set is $z_{1-\varepsilon}$ instead of $\sqrt{ \chi^2_{L; 1-\varepsilon} }$. These only coincide for $L=1$.

\section{Linearly adjustable robust counterpart: linear in what?}\label{B3}
Tractable examples of decision rules used in ARO are linear (or affine) decision rules (AARC) \citep[Chapter 14]{BenTal} or piecewise linear decision rules \citep{Chen2008}; see also Section \ref{sec:adjintro}. The AARC was introduced by \citet{BenTal:2004:AARC} as a computationally tractable method to handle adjustable variables. In the following constraint:

\[ (\vec{a} + \vec{P} \vec{\zeta}) \transp \vec{x} + \vec{b} \transp \vec{y} \leq d \qquad \forall \bm{\zeta} \in \Zset, \]
$\vec{y}$ is an adjustable variable whose value may depend on the realization of the uncertain $\vec{\zeta}$, while $\vec{b}$ does not depend on $\vec{\zeta}$ (fixed recourse). There are two different AARCs for this constraint:

 \textbf{AARC 1.} $\vec{y}$ is linear in $\vec{\zeta}$ (e.g., see \citet{BenTal:2004:AARC} and \citet[Chapter 14]{BenTal}), or

 \textbf{AARC 2.} $\vec{y}$ is linear in $\vec{a} + \vec{P} \vec{\zeta}$ (e.g., see \citet[Chapter 20.4]{AIMMS}).

\noindent Note that AARC 2 is as least as conservative as AARC 1, since the linear transformation of $\vec{\zeta} \mapsto \vec{a} +  \vec{P} \vec{\zeta}$ can only lead to loss of information, and that both methods are equivalent if the linear transformation is injective on $\Zset$. The choice for a particular method may be influenced by four factors: (i) the availability of information. An actual decision cannot depend on $\vec{\zeta}$ if $\vec{\zeta}$ has not been observed. (ii) The number of variables in the final problem. AARC 1 leads to $|\vec{\zeta}|$ extra variables compared to the RC, whereas AARC 2 leads to $|\vec{a}|$ extra variables. (iii) Simplicity for the user. Often the user observes model parameters instead of the primitive uncertainty vector. (iv) For state or analysis variables one should always use the least conservative method.

The practical issue raised in the first factor (availability of information) has been addressed with an information base matrix $\vec{P}$. Instead of being linear in $\vec{\zeta}$, $\vec{y}$ can be made linear in $\vec{P}\vec{\zeta}$. We give one example where uncertain demand is observed. Suppose there are two time periods and three possible scenarios for demand time period one and two, namely $(10,10) \transp$, $(10,11) \transp$ and $(11,11) \transp$. So, the uncertainty set of the demand vector is the convex hull of these scenarios: $\{ \vec{P} \vec{\zeta} : \vec{\zeta} \in \Zset \}$ where $\vec{P}$ is the matrix with the scenarios as columns and $\Zset = \Delta^{2} = \{\vec{\zeta} \in \R^{3} : \sum_{\ell=1}^3 \zeta_{\ell} = 1, \vec{\zeta} \geq 0\}$ is the standard simplex in $\R^3$. If the observed demand for time period one is 10, it is not possible to distinguish between $\vec{\zeta} = (1,0,0) \transp$ and $\vec{\zeta} = (0,1,0) \transp$. So, a decision for time period two can be modeled either as AARC 1 with $\vec{P} = (1,1,0)$ or as AARC 2. The latter leads to a decision rule that is easier to interpret, since it directly relates to previously observed demand.

\section{Adjustable integer variables}\label{I1}
\citet[Chapter 14]{BenTal} use parametric decision rules for adjustable continuous variables. However, their novel techniques ``generally'' cannot be applied for adjustable integer variables. In the literature two alternative approaches have been proposed. \citet{Bert2013} introduced an iterative method to treat adjustable {\it binary} variables as piecewise constant functions. The approach by \citet{Bert2010} is different and is based on splitting the uncertainty region into smaller subsets, where each subset has its own binary decision variable (see also \citet{vayanos}, and \citet{hanasu}). In this section, we briefly show this last method to treat adjustable integer variables, and show how the average behavior can be improved. We use the following notation for the general RC problem:
\begin{align*}
\text{(RC1)   }\max_{\vec{x, y, z}}\;\;  &c(\vec{x,y, z}) \\
\text{s.t.}\;\; &\vec{A(\bm{\zeta})\; x+ B(\bm{\zeta}) \;y +C(\bm{\zeta})\;z \leq d},\;\;\;\; \forall \bm{\zeta} \in \mathcal{Z},
\end{align*}
where $\vec{x} \in \R^{^{n_1}}$ and $\vec{y}\in \Z^{^{n_2}}$ are ``here and now'' variables, i.e., decisions on them are made before the uncertain parameter $\bm{\zeta}$, contained in the uncertainty set $\mathcal{Z}\subseteq\R^{^L}$, is revealed; $\vec{z}\in \Z^{^{n_3}}$ is a ``wait and see'' variable, i.e., the decision on $\vec{z}$ is made after observing (part of) the value of the uncertain parameter. $\vec{A(\bm \zeta)}\in \R^{^{m_1 \times n_1}}$ and $\vec{B(\bm{\zeta})}\in \R^{^{m_2 \times n_2}}$ are the uncertain coefficient matrices of the ``here and now'' variables. Notice that the integer ``wait and see'' variable $\vec{z}$ has an uncertain coefficient matrix $\vec{C(\bm{\zeta})} \in \R^{^{m_3 \times n_3}}$. So, unlike the ``classic'' parametric method, this approach can handle uncertainties in the coefficients of the integer ``wait and see'' variables. For the sake of simplicity, we assume that the uncertain coefficient matrices to be linear in $\bm \zeta$ and, without loss of generality, $c(\vec{x,y,z})$ is the certain linear objective function.

To model the \emph{adjustable} RC (ARC) with integer variables, we first divide the given uncertainty set $\mathcal{Z}$ into $m$ disjoint, excluding the boundaries, subsets ($\mathcal{Z}_i$, $i=1,\ldots,m$):
\begin{align*}
    \mathcal{Z} = \bigcup_{i \in \{1,\ldots,m\}} \mathcal{Z}_i,
\end{align*}
and we introduce additional integer variables $\vec{z_i}\in \Z^{^{n_3}}$ ($i=1,\ldots,m$) that model the decision in %subset
$\mathcal{Z}_i$. Then, we replicate the uncertain constraint and the objective function in (RC1) for each $\vec{z_i}$ and the uncertainty set $\Zset_i$ as follows:
\begin{align}
\text{(ARC1)   }\max_{\vec{x, y, Z}, t}\;\;  t  \nonumber\\
\text{s.t.}\;\; &c(\vec{x,y,z_i}) \geq t && \forall i \in \{1,\ldots,m\}\label{slack} \\
&\vec{A}(\bm{\zeta})\; \vec{x} + \bm{B}(\bm{\zeta})\; \bm{y} + \bm{C}(\bm{\zeta})\;\bm{z_i} \leq \bm{d} &&\forall \bm{\zeta} \in \mathcal{Z}_i, \forall i \in \{1,\ldots,m\}. \nonumber
\end{align}
Note that (ARC1) is more flexible than the non-adjustable RC (RC1) in selecting the values of integer variables, since it has a specific decision $\vec{z_i}$ for each subset $\mathcal{Z}_i$. Therefore, (ARC1) yields a robust optimal objective that is at least as good as (RC1).

Pareto inefficiency (see Section \ref{sec:solvingrc}) is also an issue in (ARC1). Hence, we must take the individual performance of the $m$ subsets into account to have a better understanding of the general performance of (ARC1). Similar to \cite{Dan}, we apply a reoptimization procedure to improve the average performance of (ARC1). More precisely, we first solve (ARC1) and find the optimal objective $t^*$. Then, we solve the following problem:
\begin{align*}
\text{(re-opt)   }\max_{\vec{x, y, Z, t}}\;\;  &\sum_{i\in\{1,\ldots,m\}} t_i \nonumber\\
\text{s.t.}\;\; &t_i \geq t^* &&\forall i \in \{1,\ldots,m\}\\
&c(\vec{x,y,z_i}) \geq t_i &&\forall i \in \{1,\ldots,m\} \\
&\vec{A(\bm{\zeta})\; x + B(\bm{\zeta})\; y + C(\bm{\zeta})\;z_i \leq d} &&\forall \bm{\zeta} \in \mathcal{Z}_i, \forall i \in \{1,\ldots,m\}, \nonumber
\end{align*}
that optimizes (i.e., maximizes) the slacks in (\ref{slack}), while the worst case objective value $t^*$ remains the same. Note that $t_i$'s are the additional variables associated with the objectives values of the subsets; (re-opt) mimics a multi-objective optimization problem that assigns equal weights to each objective, and finds Pareto efficient robust solutions.

\subsection*{Example}
Here we compare the optimal objective values of (RC1), (ARC1), and (ARC1) with (re-opt) via a toy example. For the sake of exposition, we exclude continuous variables in this example. The non-adjustable RC is given as follows:
\begin{alignat}{3}
\begin{aligned}\label{exrc}
\max_{(w, \vec{z}) \in \mathbb{Z}_+^{^3}} \;\; &5w + 3z_1 + 4z_2 \\
\text{s.t.}\;\; &(1 + \zeta_1 + 2 \zeta_2)w +(1-2\zeta_1 +\zeta_2)z_1 +(2+2\zeta_1)z_2 \leq 18 &&\forall \bm{\zeta} \in \text{Box} \\
&(\zeta_1 + \zeta_2)w +(1-2\zeta_1)z_1 +(1-2\zeta_1-\zeta_1)z_2 \leq 16 &&\forall \bm{\zeta} \in \text{Box},
\end{aligned}
\end{alignat}
where $\text{Box}=\{\zeta: -1 \leq \zeta_1 \leq 1, -1 \leq \zeta_2 \leq 1\}$ is the given uncertainty set, and $w$, $z_1$, and $z_2$ are nonnegative integer variables. In addition, we assume that $z_1$ and $z_2$ are adjustable on $\zeta_1$; i.e., the decision on these variables is made after $\zeta_1$ is being observed. Next, we divide the uncertainty set into two subsets:
\begin{align*}
&\mathcal{Z}_1 = \{(\zeta_1, \zeta_2): -1 \leq \zeta_1 \leq 0, -1 \leq \zeta_2 \leq 1\}\\
&\mathcal{Z}_2 = \{(\zeta_1, \zeta_2): 0\leq  \zeta_1 \leq 1, -1 \leq \zeta_2 \leq 1 \}.
\end{align*}
Then ARC of (\ref{exrc}) is:
\begin{align*}
\text{(Ex:ARC)   }\max_{t,w,\vec{Z}} \;\; &t \nonumber\\
\text{s.t.}\;\ &5w + 3z_1^i + 4z_2^i \geq t &&\forall i \in \{1,\ldots,m\}\\
&(1 + \zeta_1 + 2 \zeta_2)w +(1-2\zeta_1 +\zeta_2)z_1^i +(2+2\zeta_1)z_2^i \leq 18 &&\forall \bm{\zeta} \in \mathcal{Z}_i,\forall i \in \{1,\ldots,m\}\nonumber \\
&(\zeta_1 + \zeta_2)w +(1-2\zeta_1)z_1^i +(1-2\zeta_1-\zeta_1)z_2^i \leq 16 &&\forall \bm{\zeta} \in \mathcal{Z}_i, \forall i \in \{1,\ldots,m\}, \nonumber
\end{align*}
where $t\in\mathbb{R}$, $w\in\mathbb{Z}_+$, $\vec{Z}\in\mathbb{Z}^{^{2\times m}}_+$, and $m=2$ since we have two subsets. Table \ref{ex1} presents the optimal solutions of RC and ARC problems.
\begin{table}[h]\caption{RC vs ARC}\label{ex1} \centering
\begin{tabular}{c c c l} \toprule
Method &Obj. &$w$ &$\vec{z}$ \\ \midrule
RC  &29 &1 &$(z_1, z_2)=(4,3)$ \\
ARC &31 &0 &$(z^1_1, z^1_2, z^2_1, z^2_2)=(0,8,5,4)$ \\ \bottomrule
\end{tabular}
\end{table}

The numerical results show that using the adjustable reformulation we improve the objective value of the non-adjustable problem by 7\%. On the other hand, if we assume that $z_1$ and $z_2$ are adjustable on $\zeta_2$ (but not on $\zeta_1$), and we modify the uncertainty subsets $\mathcal{Z}_1$ and $\mathcal{Z}_2$ accordingly, then RC and ARC yield the same objective 29. This shows that the value of information of $\zeta_1$ is higher than that of $\zeta_2$.

Next we compare the average performance of ARC and the second stage optimization problem (re-opt) that is given by:
\begin{align*}
\max_{\vec{t},w,\vec{Z}} \; &\sum_{i\in\{1,\ldots,m\}} t_i \nonumber\\
\text{s.t.}\;\ &5w + 3z_1^i + 4z_2^i \geq t_i,\;\;\;\; t_i\geq t^* &&\forall i \in \{1,\ldots,m\} \\
&(1 + \zeta_1 + 2 \zeta_2)w +(1-2\zeta_1 +\zeta_2)z_1^i +(2+2\zeta_1)z_2^i \leq 18 &&\forall \bm{\zeta} \in \mathcal{Z}_i,\forall i \in \{1,\ldots,m\}\\
&(\zeta_1 + \zeta_2)w +(1-2\zeta_1)z_1^i +(1-2\zeta_1-\zeta_2)z_2^i \leq 16 &&\forall \bm{\zeta} \in \mathcal{Z}_i, \forall i \in \{1,\ldots,m\},
\end{align*}
where $\vec{t} \in \mathbb{R}^{^m}$. For changing the number of subsets, we again split the uncertainty sets ($\mathcal{Z}_i, i=1,\ldots,m$) on $\zeta_1$ but not on $\zeta_2$. The numerical results are presented in Table \ref{ex}.

\begin{table}[h]\caption{ARC vs re-opt for varying number of subsets}\label{ex} \centering
\begin{tabular}{c l l c c} \toprule
&\multicolumn{2}{c}{Worst Case Obj. Values per Subset} &\multicolumn{2}{c}{W.-C. Average} \\ \cmidrule(lr{.75em}){2-3} \cmidrule(lr{.75em}){4-5}
\# Subsets &\centering ARC &re-opt &ARC &re-opt\\ \midrule
1  &29 &29 &29 &29.0 \\
2  &(32, 31*)  &(34, 31*) &31.5 &32.5   \\
3  &(33, 30*, 32) &(49, 30*, 35) &31.6 &38.0\\
4  &(33, 31*, 32, 32) &(64, 34, 31*, 54) &32 &45.7 \\
5  &(33, 30*, 30*, 32, 32) &(80, 40, 30*, 33, 66) &31.4 &49.8\\
\multirow{2}{*}{8} &(32, 32, 32, 34, 31*, &(128, 64, 40, 34, 31*  &\multirow{2}{*}{32.5} &\multirow{2}{*}{61.8} \\
                   & \;\;33, 33, 33) & \; 36, 54, 108) & & \\
\multirow{2}{*}{10} &(32, 32, 32, 32, 34, &(160, 80, 52, 40, 34,  &\multirow{2}{*}{32.5} &\multirow{2}{*}{64.3} \\
                    & \;\;31*, 33, 33, 33, 33) & \;\;31*, 33, 45, 66, 135) & &  \\ \bottomrule
\multicolumn{5}{l}{($\ast$) denotes the worst case (w.-c.) objective value over all subsets}
\end{tabular}
\end{table}

The first column of the table presents the number of subsets used in ARC, and we assume that the domain of $\zeta_1$ is divided into equally sized intervals (e.g., if the number of subsets is three, then the intervals are $-1 \leq \zeta_1 \leq -0.33, -0.33 \leq \zeta_1 \leq 0.33$, and $ 0.33 \leq \zeta_1 \leq 1$). The second column reports objective values of the subproblems associated with the subsets in ARC. The third column presents the objective values of the subproblems when we apply (re-opt). The fourth and fifth columns show the averages of the results in columns two and three. As anticipated, intra row comparisons show that ARC and (re-opt) yield the same worst case performance for a fixed number of subsets, and (re-opt) yields significantly better average performances than ARC. Moreover, the average performance improves with the number of subsets. Notice that the average performance of the RC solution is not reported in Table \ref{ex} because it has the same average performance, that is 29, for any given number of subsets. Nevertheless, it is important to see the significant average objective value improvement made by ARC with (re-opt) for the ``fixed'' performance of the RC. Last but not least, the optimal objective value 31, which is obtained when the number of subsets is two, four, eight and ten in Table \ref{ex}, is the global optimal of the ARC; for details on optimality see the following section where this example will be revisited.

\subsection*{Optimality}
To quantify how far is the optimal objective value ($t^*$) of (ARC1) from that of the best possible solution, we need to define an efficient lower bound (or an upper bound for a maximization problem) for the best objective. One way of finding such a bound is by solving (RC1), where the uncertainty set $\Zset$ is replaced with a finite subset (denoted by $\hat{\Zset}$), and where each scenario from $\hat{\Zset}$ has a separate second-stage decision \citep{hadji,Bert2013,postek}. The optimal objective value of such a formulation is always a lower bound for the best possible objective value, since it is an optimal (i.e., not restricted to an affine form) adjustable solution for a smaller uncertainty set. More precisely, the lower bound problem is given as follows:
\begin{align*}
\text{(BRC)    }\min_{\vec{x, y, z^{(\bm{\zeta})}}, t^{\text{lb}}}\;\;  t^{\text{lb}}  \\
\text{s.t.}\;\; &c(\vec{x,y,z^{(\bm{\zeta})}}) \leq t^{\text{lb}} && \forall \bm{\zeta} \in \hat{\Zset} \\
&\vec{A(\bm{\zeta})\; x + B(\bm{\zeta})\; y +C(\bm{\zeta})\;z^{(\bm{\zeta})} \leq d} &&\forall \bm{\zeta} \in \hat{\Zset}
\end{align*}
where $\hat{\Zset}$ is a finite subset of $\Zset$, as explained above. Now the question that has to be answered is: how to construct $\hat{\Zset}$ efficiently? \citet{postek} proposes to first find the optimal solution of (ARC1) for a given number of subsets, and then taking $\hat{\Zset}$ as uncertain parameters that maximizes the left-hand side of at least one constraint. For additional details on improving the lower bound we refer to \citet[Chapter~4.2]{postek}.

\emph{Example (Ex:ARC) revisited.} The solution of (Ex:ARC) for two subsets (i.e., $m=2$) is given in the second row of Table \ref{ex1}. The associated finite ``worst case'' subset for this solution is $\hat{\Zset}=\{(0,1),(0,-1)\}$, and the upper bound for the best possible worst case objective is $t^{\text{ub}}=31$ (this is obtained by solving the upper bound reformulation of (BRC) for $\hat{\Zset}$). Therefore, the optimal objective value of (Ex:ARC) is bounded above by 31 for any given number of subsets. So, the solution found for two subsets (see Table \ref{ex}) is optimal w.r.t. the worst case.

\subsection*{Tractability}
It is important to point out that our adjustable reformulation and the ``non-adjustable'' RC have the same ``general'' mathematical complexity, but the adjustable reformulation increases the number of variables and constraints by a factor $m$ (the number of subsets), so that if the number of integer variables is high (say a few hundreds) then the resulting adjustable RC may be intractable. Dividing the main uncertainty set $\mathcal{Z}$ into more subsets $\mathcal{Z}_i$ may improve the objective value by giving more freedom in making adjustable decisions, but the decision maker should make the tradeoff between optimality and computational complexity.

\section{Robust counterparts of equivalent deterministic problems are not necessarily equivalent}\label{D2}

In this section we show that the robust counterparts of equivalent deterministic problems are not always equivalent.
The message in this section is thus that one has to be careful with reformulating optimization problems,
since the corresponding robust counterparts may not be the same.

\medskip\noindent
Let us start with a few simple examples.
The first one is similar to the example in \citet[p.~13]{BenTal}.
Consider the following constraint:
\[  (2+\zeta)x_1 \leq 1,  \]
where $\zeta$ is an (uncertain) parameter.
This constraint is equivalent to:
\begin{equation*}
\begin{cases}
  (2+\zeta)x_1 + s = 1  \\
  s \geq 0.
\end{cases}
\end{equation*}
However, the robust counterparts of these two constraint formulations, i.e.
\begin{equation}  (2+\zeta)x_1 \leq 1 \quad \forall \zeta \ : \ | \zeta | \leq 1,  \label{RC} \end{equation}
and
\begin{equation}
\begin{cases}
  (2+\zeta)x_1 + s = 1 \quad \forall  \zeta \ : \  |\zeta | \leq 1\\
  s \geq 0,
\end{cases} \label{RC2}
\end{equation}
in which the uncertainty set for $\zeta$ is the set $\{\zeta: |\zeta| \leq 1\}$, are not equivalent.
It can easily be verified that the feasible set for robust constraint (\ref{RC}) is:  $x_1 \leq 1/3$,
while for the robust constraint (\ref{RC2}) this is $ x_1 =0$.
The reason why (\ref{RC}) and (\ref{RC2}) are not equivalent is that by adding the slack variable, the inequality becomes an equality that has to be satisfied for all values of the uncertain parameter, which is very restrictive. The general message is therefore: {\it do not introduce slack variables  in uncertain constraints, unless they are adjustable like in \citet{Daniel},  and avoid uncertain equalities}.

\medskip\noindent
Another example is the following constraint:
\begin{equation*}  |  x_1 - \zeta | + |  x_2 -\zeta | \leq 2, \end{equation*}
which is equivalent to:
\begin{equation*}
  \begin{cases}
   y_1 + y_2 \leq 2  \\
   y_1 \geq  x_1 - \zeta   \\
   y_1 \geq  \zeta  -x_1  \\
   y_2 \geq x_2 -\zeta \\
   y_2 \geq  \zeta - x_2.
     \end{cases}
\end{equation*}
However, the robust versions of these two formulations, namely:
\begin{equation}  |  x_1 - \zeta | + |  x_2 -\zeta | \leq 2 \quad \forall \ \zeta \ : \ | \zeta | \leq 1,    \label{Ex2}\end{equation}
and:
\begin{equation}
  \begin{cases}
   y_1 + y_2 \leq 2     \\
   y_1 \geq  x_1 - \zeta   \quad \forall \zeta \ : \ | \zeta | \leq 1\\
   y_1 \geq  \zeta  -x_1   \quad \forall \zeta \ : \ | \zeta | \leq 1 \\

   y_2 \geq x_2 -\zeta    \quad\forall \zeta \ : \ | \zeta | \leq 1\\
   y_2 \geq  \zeta - x_2   \quad\forall \zeta \ : \ | \zeta | \leq 1,
     \end{cases}
   \label{ExR2}
\end{equation}
are not equivalent. Indeed, it can easily be checked that the
set of feasible solutions for  (\ref{Ex2}) is $(\theta,-\theta)$, $-1 \leq  \theta \leq 1$, but the
only feasible solution for (\ref{ExR2}) is $\vec{x} = (0,0)$. The reason for this is that in (\ref{ExR2}) the uncertainty is split over several constraints, and since the concept of RO is constraint-wise, this leads to different problems, and thus different solutions. The following linear optimization reformulation, however, is equivalent to (\ref{Ex2}):
\begin{equation}
  \begin{cases}
  x_1 - \zeta   +   x_2 -\zeta  \leq 2 &\quad \forall \ \zeta \ : \ | \zeta | \leq 1  \\
  x_1 - \zeta   + \zeta - x_2 \leq 2 &\quad \forall \ \zeta \ : \ | \zeta | \leq 1   \\
 \zeta -x_1+   x_2 -\zeta  \leq 2 &\quad \forall \ \zeta \ : \ | \zeta | \leq 1  \\
 \zeta - x_1  +\zeta - x_2  \leq 2 &\quad \forall \ \zeta \ : \ | \zeta | \leq 1.
\end{cases}
\end{equation}

The general rule therefore is: {\it do not split the uncertainty in one constraint over more constraints, unless the uncertainty is disjoint.}
In particular do not use ``definition variables'' if this leads to such a splitting of the uncertainty.

\medskip\noindent
In the remainder we give a general treatment of some often used reformulation tricks to reformulate nonlinear problems into linear ones, and discuss whether the robust counterparts are equivalent or not.

\medskip\noindent
\begin{itemize}
\item {\bf Maximum function.}
Consider the following constraint:
\[ \vec{a(\bm{\zeta})\transp  x}  +  \max_k \vec{b_k(\bm{\zeta})\transp  x} \leq d(\bm{\zeta}) \quad \forall \bm{\zeta} \in \Zset, \]
where $\bm{\zeta} \in \Zset$ is the uncertain parameter, and $\vec{a(\bm{\zeta})}$, $\vec{b_k(\bm{\zeta})}$, and $d(\bm{\zeta})$ are parameters that depend linearly on $\bm{\zeta}$.
A more conservative reformulation for this constraint is:
\[ \left\{ \begin{array}{l l}
       \vec{a(\bm{\zeta})\transp  x } +  z  \leq d(\bm{\zeta}) \quad \quad &\forall \bm{\zeta} \in \Zset  \\
       z  \geq  \vec{b_k(\bm{\zeta})\transp  x} &\forall k, \forall \bm{\zeta} \in \Zset,
   \end{array} \right. \]
since the uncertainty is split over more constraints.
The exact reformulation is:
\[ \vec{a(\bm{\zeta})\transp  x  +  b_k(\bm{\zeta})\transp  x} \leq d(\bm{\zeta}) \quad \forall k, \forall \bm{\zeta} \in \Zset. \]
Note that in many cases we have  ``a sum of max'':
\[ \vec{a(\bm{\zeta})\transp  x}  +  \sum_i \max_k \vec{b_{ik}(\bm{\zeta}) \transp x} \leq d(\bm{\zeta}) \quad \forall \bm{\zeta} \in \Zset. \]
Important examples that contain such constraints are production-inventory problems.
We refer to \citet{sumofmax} for an elaborate treatment on exact and approximate reformulations of such constraints.

\medskip\noindent
\item {\bf Absolute value function.}
Note that $|x| = \max\{x,-x\}$, and hence this is a special case of the $\max$ function, treated above.

\medskip\noindent
\item {\bf Linear fractional  program.}
Consider the following robust linear fractional problem:
\begin{equation} \left\{ \begin{array}{l l}
\displaystyle \min_{\vec{x}} &\max_{\bm{\zeta} \in \Zset}\frac{\alpha(\bm{\zeta}) + \vec{c(\bm{\zeta})\transp  x}}{\beta(\bm{\zeta}) + \vec{d(\bm{\zeta})\transp  x}} \\
\text{s.t.} &\displaystyle \sum_j a_{ij}x_j  \geq b_i \quad \forall i  \\
&\vec{x \geq 0}, \end{array} \right.  \label{RLF}  \end{equation}
where $\bm{\alpha}(\bm{\zeta})$, $\vec{c}(\bm{\zeta})$, $\bm{\beta}(\bm{\zeta})$, and $\vec{d}(\bm{\zeta})$ are parameters that depend linearly on $\bm{\zeta}$. Moreover, we assume that $\beta(\bm{\zeta}) + \vec{d}(\bm{\zeta})\transp  \vec{x} > 0$, for all feasible $\vec{x}$ and for all $\bm{\zeta} \in \Zset$.
For the non-robust version one can use the Charnes-Cooper transformation that is proposed by \citet{frac} to obtain an equivalent linear optimization problem. However, if we apply this transformation to the robust version, we obtain:
\[ \left\{ \begin{array}{r l l}
\displaystyle \min_{\vec{y},t} &\displaystyle\max_{\bm{\zeta} \in \Zset}  \alpha(\bm{\zeta})t + \vec{c(\bm{\zeta})\transp  y} \\
\text{s.t.} &\beta(\bm{\zeta})t + \vec{d(\bm{\zeta})\transp  y} = 1 &\quad \forall \bm{\zeta} \in \Zset  \\
&\displaystyle \sum_j a_{ij}y_j  \geq b_i t &\quad \forall i\\
&\vec{y \geq 0}, \ t\geq 0,  \end{array}\right.\]
which is not equivalent to (\ref{RLF}) since the uncertainty in the original objective is now split over the objective and a constraint.
A better way to deal with such problems is to solve the robust linear problem
\[\left\{ \begin{array}{r l}
\displaystyle \min_{\vec{x}} &\displaystyle \max_{\bm{\zeta} \in \Zset} \left[\alpha(\bm{\zeta}) + \vec{c(\bm{\zeta})\transp  x}  - \lambda \left( \beta(\bm{\zeta}) + \vec{d(\bm{\zeta})\transp  x} \right)  \right] \\
\text{s.t.} &\displaystyle \sum_j a_{ij}x_j  \geq b_i \\
&\vec{x \geq 0},  \end{array} \right. \]
for a fixed value of $\lambda$, and then find the minimal value of $\lambda$ for which this optimization problem still has a  non positive optimal value. One can use for example binary search on $\lambda$ to do this. For a more detailed treatment of robust fractional problems we refer to \citet{BramFrac}.

\medskip\noindent
\item {\bf Product of binary variables.}
Suppose that a robust constraint contains a product of binary variables, say $xy$, with $x,y \in \{0,1\}$. Then one can use the standard way to linearize this:
\[ \left\{ \begin{array}{l}
       z \leq x \\
       z \leq y  \\
       z \geq x + y -1 \\
       z\geq 0,
   \end{array} \right. \]
and replace $xy$ with $z$.
One can use this reformulation since the added constraints do not contain uncertain parameters.

\medskip\noindent
\item {\bf Product of binary and continuous variable.}
A product of a binary and a continuous variable that occurs in a robust constraint can also be reformulated in linear constraints, in a similar way as above.
However, note that in the following robust constraint:
\[ \vec{a(\bm{\zeta})\transp  x} + z \vec{b(\bm{\zeta})\transp  x}  \leq d(\bm{\zeta}) \quad \forall \bm{\zeta} \in \Zset, \]
where $z\in \{0,1\}$,
one cannot use the standard trick:
\begin{equation} \left\{ \begin{array}{l l}
       \vec{a(\bm{\zeta})\transp  x} + z y  \leq d(\bm{\zeta})  &\quad\forall \bm{\zeta} \in \Zset \\
           y \geq  \vec{b(\bm{\zeta})\transp  x}   &\quad \forall \bm{\zeta} \in \Zset,
\end{array} \right. \label{zy} \end{equation}
and then linearize $zy$. This is not possible since in (\ref{zy}) the uncertainty is split over different constraints.
A correct reformulation is:
\begin{equation} \left\{ \begin{array}{l l}
       \vec{a(\bm{\zeta})\transp  x} +   \vec{b(\bm{\zeta})\transp  x} \leq d(\bm{\zeta}) + M (1-z)  &\quad\bm{\zeta} \in \Zset \\
        \vec{a(\bm{\zeta})\transp  x} \leq d(\bm{\zeta})+Mz  &\quad\bm{\zeta} \in \Zset,
   \end{array} \right. \end{equation}
where $M$ is a sufficiently big number.
\medskip\noindent
\item {\bf $K$ out of $N$ constraints should be satisfied.}
Suppose the restriction is that at least $K$ out of the $N$ robust constraints
\begin{equation}
\vec{a_i(\bm{\zeta})\transp  x} \leq d_i(\bm{\zeta}) \quad \forall \bm{\zeta} \in \Zset
\label{eq:22}
\end{equation}
should be satisfied, where $i\in\{1,\ldots,N\}$. Then one can use the standard way
\[ \left\{ \begin{array}{l l}
       \vec{a_i(\bm{\zeta})\transp  x }\leq d_i(\bm{\zeta})  + M (1-z_i) &\quad \forall \bm{\zeta} \in \Zset,  \forall i \\
       \displaystyle \sum_i z_i \geq K \\
       z_i  \in  \{0,1\}  &\quad \forall i,
\end{array} \right. \]
where $M$ is a sufficiently big number.
However, if the restriction is that $\forall \bm{\zeta} \in \Zset$ at least $K$ out of the $N$ constraints should be satisfied (notice the difference with (\ref{eq:22})), then the above constraint-wise formulation is not equivalent and is overly conservative.
We do not see how to model such a constraint correctly. Maybe an adversarial approach could be used for such constraints.

\medskip\noindent
\item {\bf If-then constraint.}
Since an ``if-then constraint'' can be modeled as an at least $1$ out of $2$ constraints, the above remarks hold.
\end{itemize}

\medskip\noindent
Up to now we only described linear optimization examples. Similar examples can be given for conic and nonlinear optimization. In \citet{lobo} for example, many optimization problems are given that can be modeled as conic quadratic programming problems. However, for many of them it holds that the corresponding robust counterparts are not the same. This means that if an optimization problem is conic quadratic representable, then the robust counterparts are not automatically the same, and hence in such cases the robust optimization techniques for CQP cannot be used.

\section{How to deal with equality constraints?}\label{D1}

Equality constraints containing uncertain parameters should be avoided as much as possible, since often such constraints restrict the feasibility region drastically or even lead to infeasibility. Therefore, the advice is:  {\it do not use slack variables unless they are adjustable, since using slack variables leads to equality constraints}; see \citet[Chapter 2]{BenTal}. However, equality constraints containing uncertain parameters cannot always be avoided.
There are several ways to deal with such uncertain equality constraints:
\begin{itemize}
\item In some cases it might be possible to convert the equality constraints into inequality constraints. An illustrating example is the transportation problem: the demand constraints can either be formulated as equality constraints or as inequality constraints. The structure of the problem is such that at optimality these inequalities are tight.
\item  The equality constraints can be used to eliminate variables. This idea is mentioned in \citet{BenTal}. However, several questions arise. First of all, after elimination of variables and after the resulting problem has been solved, it is unclear which values to take for the eliminated variables, since they also depend on the uncertain parameters. This is no problem if the eliminated variables are  adjustable variables or state/analysis variables, since there is no need to know their optimal values. A good example  is the production-inventory problem for which one can easily eliminate the state or analysis variables indicating the inventory in different time periods. See, e.g., \citet{BenTal}. Secondly, suppose the coefficients with respect to the variables that will be eliminated contain uncertain parameters. Eliminating such variables leads to problems that contain non-linear uncertainty, which are much more difficult to solve.
To illustrate this, let us consider the following two constraints of an optimization problem:
\[  \zeta_1 x_1 + x_2 + x_3 = 1, \quad  x_1 + x_2  + \zeta_2x_3  \leq 5,\]
in which $\zeta_1$ and $\zeta_2$ are uncertain parameters. Suppose that $x_1$, $x_2$ and $x_3$ are all adjustable in $\zeta_1$. Then there are three options for elimination:
\begin{enumerate}
      \item {\it Elimination of $x_1$.} Let us assume that $\zeta_1=0$ is not in the uncertainty set. By substituting $x_1 = (1-x_2-x_3)/\zeta_1$ the inequality becomes:
  \[ \left( 1- \frac{1}{\zeta_1}\right) x_2 + \left(\zeta_2- \frac{1}{\zeta_1}\right) x_3  \leq 5- \frac{1}{\zeta_1}. \]
The disadvantage of eliminating $x_1$ is thus that the uncertainty in the inequality  becomes nonlinear.
   \item {\it Elimination of $x_2$.} By substituting $x_2 = 1-\zeta_1x_1-x_3$ the inequality becomes:
\[ (1-\zeta_1)x_1 + (\zeta_2-1)x_3  \leq 4, \]
which is linear in the uncertain parameters.
\item {\it Elimination of $x_3$}. By substituting $x_3 = 1-\zeta_1x_1 - x_2$ the inequality becomes:
\[ (1-\zeta_1\zeta_2)x_1  +  (1-\zeta_2)x_2   \leq 5-\zeta_2, \]
which is nonlinear in the uncertain parameters.
We conclude that from a computational point of view it is more attractive to eliminate $x_2$.
\end{enumerate}
It is important to note that different choices of variables to eliminate may lead to different optimization problems.

\item  If the constraint contains state or analysis variables one could make these variables adjustable and use decision rules, thereby introducing much more flexibility. One can easily prove that when the coefficients for such variables in the equality constraint do not contain uncertain parameters and the equality constraint is linear in the uncertain parameters, then using linear decision rules for such variables is equivalent to eliminating these variables. To be more precise: suppose the linear equality constraint is
\[ \vec{q(\bm{\zeta})\transp  x} + y  = r, \]
where $\vec{q(\bm{\zeta})}$ is linear in $\bm{\zeta}$, and $y$ is an state or analysis variable (without loss of generality we assume the coefficient for $y$ is $1$). Then it can easily be proven that substituting $y=r - \vec{q(\bm{\zeta})\transp  x}$ everywhere in the problem is equivalent to using a linear decision rule for $y$. To reduce the number of extra variables, it is therefore better to eliminate such variables.
\end{itemize}

\section{On maximin and minimax formulations of RC}\label{I2}
In this section, we consider an uncertain LP of the following general form:
\begin{equation*}
\text{(LP) } \max_{\vec{x}\geq \vec{0}}\; \{\vec{c\transp   x} : \vec{Ax} \leq \vec{d}\},
\end{equation*}
where without loss of generality $\vec{A}$ is the uncertain coefficient matrix that resides in the uncertainty set $\Uset$. So the general RC is given by
\begin{equation*}
\text{(R-LP) } \max_{\vec{x}\geq \vec{0}}\; \{\vec{c\transp   x} : \vec{Ax} \leq \vec{d}\;\; \forall \vec{A} \in \Uset\}.
\end{equation*}
Here we show that (R-LP) can be reformulated as:
\begin{equation*}
\text{(RF) }\min_{\vec{A} \in \Uset} \max_{\vec{x}\geq \vec{0}} \; \{\vec{c\transp   x} : \vec{Ax} \leq \vec{d}\},
\end{equation*}
if the uncertainty is constraint-wise; however if this condition is not met, then (RF) may not be equivalent to (R-LP).
\begin{remark}
This shows that the statement ``RO optimizes for the worst case $\vec{A}$'' is too vague. Also the maximin reformulation:
\begin{equation*}
\max_{\vec{x}\geq \vec{0}} \min_{\vec{A} \in \Uset}  \; \{\vec{c\transp   x} : \vec{Ax} \leq \vec{d}\},
\end{equation*}
is usually not equivalent to (R-LP). This is because we can almost always find an $\vec{x}\geq 0$ such that no $\vec{A}\in \Uset$ exists for which $\vec{Ax \leq d}$; therefore, we minimize over an empty set, and have $+\infty$ for the maximin objective. Also when $\vec{x}$ is selected such that at least one feasible $\vec{A}$ exists (e.g., see \citet{falk}), it is easy to find examples where both formulations are not equivalent.
\end{remark}
To show (R-LP)=(RF) when the uncertainty is constraint-wise, we first take the dual of the (inside) maximization problem of (RF) [$\max_{\vec{x}\geq \vec{0}} \vec{c\transp   x} : \vec{Ax} \leq \vec{d}$]. Then, substituting the dual with the primal (maximization) problem in (RF) gives:
\begin{equation*}
\text{(OC-LP) }\min_{\vec{A} \in \Uset, \vec{y}\geq \vec{0}}\; \{\vec{d\transp y} : \vec{A^T y} \geq \vec{c}\},
\end{equation*}
where val(OC-LP) = val(RF) holds under regularity conditions at optimality. Note that the constraints of (RF) can be formulated as [$\vec{a_i^{T}x} \leq d_i, \forall \vec{a_i} \in \Uset_i, i =1,\ldots,m$], if the uncertainty is constraint-wise. \citet[Theorem 3.1]{Beck} and \citet{soyster2013unifying} show that (OC-LP)---which is the \emph{optimistic counterpart of the dual problem}---is equivalent to the dual of the general robust counterpart (R-LP), and val(OC-LP) = val(R-LP) holds  for constraint-wise uncertainty, and disjoint $\Uset_i$'s. Therefore, val(R-LP) = val(RF) is also satisfied when the uncertainty is constraint-wise. However, if (some of) the uncertainty at (some of) the constraints are dependent in (R-LP), then we may not sustain the associated equivalence. The following example shows such a situation.

\subsection*{Example}
Consider the following toy RC example in which the uncertainty is not constraint-wise:
\begin{align*}
\text{(RC-Toy) }\max_{\vec{y}} \;\; &y_1 + y_2 \\
\text{s.t. }      &a_1 y_1 \leq 1, \ \ a_2 y_2 \leq 1 \;\; \forall \vec{a} \in \R^{^2}:\twonorm{\vec{a}} \leq 1,
\end{align*}
where two constraints of the problem are dependent on each other via the ellipsoidal uncertainty set [$\vec{a} \in \R^{^2}:\twonorm{\vec{a}} \leq 1$].
The robust reformulation of the (RC-Toy) is as follows:
\begin{align*}
\text{(RF-Toy) }\min_{\vec{a}:\twonorm{\vec{a}} \leq 1} \max_{\vec{y}}  \;\; &y_1 + y_2 \\
\text{s.t. }      &a_1 y_1 \leq 1, \ \ a_2 y_2 \leq 1,
\end{align*}
and the optimistic counterpart (OC) of the problem is
\begin{align*}
\text{(OC-Toy) }\min_{\vec{x\geq 0,\; a:\twonorm{\vec{a}}\leq 1}} \;\; &x_1 + x_2 \\
\text{s.t. }      &a_1 x_1 = 1, \ \ a_2 x_2 = 1.
\end{align*}
(RC-Toy) attains an optimal objective value of 2, whereas the (RF-Toy)'s optimal objective value is $2\sqrt{2}$.
Therefore, the robust reformulation (RF-Toy) is not equivalent to the general RC problem (RC-Toy) in this situation.
However, val(RF-Toy) = val(OC-Toy) from duality.

\section{Quality of robust solution}\label{B2}
In this section we describe how to assess the quality with respect to robustness of a solution based on a simulation study. We first identify four focus points for performing a Monte Carlo experiment, and conclude with two statistical tests that can be used to compare two solutions.\\
\\
{\bf Choice of the uncertainty set.}
For a comparison between different solutions, it is necessary to define an uncertainty set $\Uset$ that is used for evaluation. This set should reflect the real-life situation. The uncertainty set that is used for optimization may be different than the set for evaluation. For example, an ellipsoidal set may be used to reduce the conservatism when the real-life uncertainty is a box, while still maintaining a large probability of constraint satisfaction \cite[p.~34]{BenTal}. \\
\\
{\bf Choice of the probability distribution.}
A simulation requires knowledge of the probability distribution on the uncertainty set. If this knowledge is ambiguous, it may be necessary to verify whether the simulation results are sensitive with respect to changes in this distribution. For example, \citet{Rozenblit} performs different simulations, each based on a probability distribution with a different skewness level. \\
\\
{\bf Choice of the sampling method.}
For univariate random variables it is computationally easy to draw a random sample from any given distribution. For multivariate random variables rejection sampling can be used, but it may be inefficient depending on the shape of the uncertainty set, e.g., for an uncertainty set with no volume. A more efficient method for sampling from an arbitrary continuous probability distribution is ``hit and run'' sampling \citep{Belisle1993}. An R package for uniform hit and run sampling from a convex body is also available \citep{hitandrun}. \\
\\
{\bf Choice of the performance characteristics.}
From a mathematical point of view there is no difference between uncertainty in the objective and uncertainty in the constraints since an uncertain objective can always be reformulated as a certain objective and an uncertain constraint. However, the distinction between an uncertain objective and an uncertain constraint is important for the interpretation of a solution. First, we look at the effects of adjustable RO and reformulations, then we present the performance characteristics. \\
\\
{\it Effect of adjustable RO.}
When one or more ``wait and see'' variables are modeled as adjustable variables, uncertain parameters may enter the objective function. In that case the performance characteristics for uncertainty in the objective become applicable. \\
\\
{\it Effect of reformulations.}
Reformulations are sometimes necessary to end up with a tractable model. The evaluation should be based on the original model, since reformulations introduce additional constraints whose violation is not necessarily a problem. Take for example an inventory model that has constraints on variables that indicate the cost at a certain time period (e.g., constraints \eqref{inv:1} and \eqref{inv:2}). These constraints have been introduced to model the costs in the objective function. A violation of these constraints does not render the solution infeasible but does affect the objective value (i.e., the costs of carrying out the solution). \\
\\
{\it Performance characteristics for uncertainty in the constraints.}
For an uncertain constraint $f(\vec{a}, \vec{\zeta}) \leq 0$ for all $\vec{\zeta}$ in $\Zset$, the violation is $\max\{0, f(\vec{a}, \vec{\zeta}) \}$. Meaningful statistics are the probability on positive violation and the distribution of the violation (average, worst case, standard deviation) under the condition that the violation is positive. When multiple constraints are uncertain, these statistics can be computed per constraint. Additionally, the average number of violated constraints can be reported.

There is a clear trade-off between the objective value and constraint violations. The difference between the worst case objective value of the robust solution and the nominal objective value of the nominal solution is called the \emph{price of robustness} (PoR) \citep{BertsimasPoR}. It is useful if the objective is certain, since in that case PoR is the amount that has to be paid for being robust against constraint violations. We observe that PoR is also used when the objective is uncertain. We discourage this, since it compares the nominal solution in case there is no uncertainty with the robust solution where the worst case occurs, so it compares two different scenarios. \\
\\
{\it Performance characteristics for uncertainty in the objective.}
Uncertainty in the objective affects the performance of a solution. For every simulated uncertainty vector, the actual objective value can be computed. One may be interested in the worst case, but also in the average value or the standard deviation. For a solution that is carried out many times, reporting the average performance is justified by the law of large numbers. The worst case may be more relevant when a solution is carried out only once or a few times, e.g., when optimizing a medical treatment plan for a single patient. These numbers show what objective value to expect, but they do not provide enough information about the quality of a solution since a high standard deviation is not necessarily undesirable. A robust solution is good when it is close to the \emph{perfect hindsight} (PH) solution. The PH solution is the solution that is obtained by optimizing the decision variables for a specific uncertainty vector as if it is fully known beforehand. This has to be done for every simulated uncertainty vector, and yields an utopia solution. The PH solution may have a large variation, causing a high variation of good solutions as well. \\
\\
{\it Performance characteristics for any problem.}
Regardless of whether the uncertainty is in the objective or in the constraints, the mean and associated standard deviation of the difference between the actual performance of a solution and the PH solution are useful for quantifying the quality of a solution. The mean difference between the PH solution and a fully robust solution is defined as the \emph{price of uncertainty} (PoU) by \citet{BenTal:2005:retailer}. It is the maximum amount that a company should invest for reducing the level of uncertainty, e.g., by using more accurate forecasting techniques. It can also be interpreted as the regret of choosing a certain solution rather than the PH solution. Alternative names for PoU are ``cost of robustness'' \citep{Gregory2011} or ``price of robustness'' \citep{BenTal:2004:AARC}, which are less descriptive than ``price of uncertainty'' and may cause confusion with price of robustness from \citep{BertsimasPoR}. A low mean PoU and a low standard deviation characterize a good solution.

Subtracting the mean objective value of the nominal solution from the mean value of a robust solution yields the \emph{actual price of robustness} (APoR) \citep{Rozenblit}. APoR can be interpreted as the expected price that has to be paid for using the robust solution rather than the nominal solution, which is negative if RO offers a solution that is better on average. PoR equals APoR when uncertainty only occurs in the constraints.

For multistage problems one may also follow a folding horizon (FH) approach. With FH in each stage where a part of the uncertain parameter is observed, that information is used to optimize for the remaining time periods. This is done by taking the original optimization problem, fixing the decision variables for previous stages, and fixing the elements of the uncertain parameter that have been observed. This allows a fair comparison between a dynamic solution (e.g., created by the AARC) and a static solution (e.g., the nominal solution) when in real-life the static solution is reoptimized in every stage.\\
\\
{\bf Comparing two solutions.}
We provide several comparison criteria and provide the corresponding statistical test to verify whether one solution is better than another solution. The tests will be demonstrated in Section \ref{B4}. We will assume that the data for the statistics test is available as $n$ pairs $(X_i,Y_i)$ ($i=1,2,\ldots,n$), where $X_i$ and $Y_i$ are performance characteristics in the $i$'th simulation. For uncertainty in the objective, they can be objective values whereas for uncertainty in the constraints they can be the numbers of constraint violations or the sizes of the constraint violations. We assume that $(X_i,Y_i)$ and $(X_j,Y_j)$ are independent if $i \neq j$, and that smaller values are better. When a conjecture for a test is based on the outcome of a simulation study, the statistical test must be performed with newly generated data to avoid statistical bias. While for the statistical tests it is not necessary that $X_i$ and $Y_i$ are based on the same simulated uncertainty vector $\vec{\zeta}$, it increases the power of the test since $X_i$ and $Y_i$ will be positively correlated. This reduces the variance of the difference: $\Var(X_i - Y_i) = \Var(X_i) + \Var(Y_i) - 2 \Cov(X_i,Y_i)$, which is used in the following tests:
\begin{itemize}
\item The sign test for the median validates the null hypothesis that the medians of the distributions of $X_i$ and $Y_i$ are equal. This tests the conjecture that the probability that one solution outperforms the other is larger than 0.5.
\item The t-test for the mean validates the null hypothesis that the means of the distributions of $X_i$ and $Y_i$ are equal. This tests the conjecture that one solution outperforms the other in long run average behavior.
\end{itemize}

\section{RC may take better ``here and now'' decisions than AARC}\label{B4}
A linear decision rule is a linear approximation of a more complicated decision rule. It dictates what to do at each stage as a linear function of observed uncertain parameters, but it is not guaranteed to be the optimal strategy. Every time a decision has to be made it is possible to either follow the linear decision rule, or to reoptimize the AARC for the remaining time periods based on everything that is observed up till then. We will refer to the latter as the AARC-FH, where FH stands for folding horizon. \citet{BenTal:2005:retailer} compare the AARC with the AARC-FH, and show that the latter produces better solutions on average. A comparison that involves AARC-FH assumes that there is time to reoptimize. It is therefore natural to also make a comparison with the RC-FH, where the RC is solved for the full time horizon and re-optimized for the remaining time period every time a part of the uncertain parameters is unveiled. On average, the RC-FH may outperform the AARC \citep{Cohen2007,Rozenblit}.

In the remainder of this section we will evaluate both the average and the worst case performance of the nominal solution with FH, the RC-FH and the AARC-FH. A comparison between RC-FH and AARC-FH is new, and shows which model takes the best ``here and now'' decisions.

We first give an example for the worst case performance. Consider a warehouse that transfers one good. The current inventory is $x_0 = 5$, the holding costs per time period are $h = 1$, the backlogging costs per time period are $b = 2$. In the first period, any nonnegative (not necessarily integer) amount can be ordered while in the second period the maximum order quantity is $q_2^{max}=3$. Let $T = \{1,2\}$, let $q_t$ be the order quantity in time period $t$, and let $c_t$ denote the costs associated with time period $t$. The ending inventory can be returned to the supplier without penalty fee at time period three. The optimization problem can be formulated as:

\begin{align}
& \min  			\qquad	&& \sum_{t \in T} c_t &&                                           \notag \\
& \mbox{s.t.} \qquad	&& c_t \geq  (x_0 + \sum_{i=1}^T q_i - d_i)h && \forall t \in T  \label{inv:1}\\
&                     && c_t \geq -(x_0 + \sum_{i=1}^T q_i - d_i)b && \forall t \in T  \label{inv:2} \\
&                     && q_2 \leq q_2^{max}                                              \notag \\
&                     && q_t \in \R_+                                && \forall t \in T. \notag
\end{align}
Suppose the demand vector $\vec{d}$ is uncertain but is known to reside in a ball around $\vec{5}$ with radius $5$. We will use this uncertainty set both for optimization and for evaluation.

\begin{figure}
	\centering
	\caption{Approximation of the total worst case costs for an FH strategy as a function of the initial order quantity $q_1$.}
		\includegraphics[scale=1]{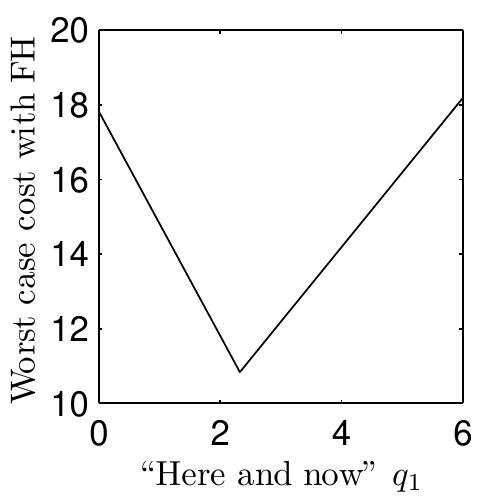}
    \label{fig:inv-worstcase}
\end{figure}

For this small example, it is possible to approximate the worst case costs for an FH approach as a function of the ``here and now'' decision $q_1$ as follows. For each $q_1$ in a range of values, we have randomly drawn $100$ uncertain demand vectors from the boundary of the uncertainty set. For each demand vector we have computed the inventory level at the beginning of the second time period ($=x_0 + q_1 - d_1$). Based on this inventory level,  we reoptimized the order quantity for the second time period, where $d_2$ was assumed to reside in the interval $[5-r, 5+r]$ with $r=\sqrt{25 - (d_1-5)^2}$ (so that the full $\vec{d}$ vector is in a ball around $\vec{5}$ with radius $5$). Then we computed the total costs over both time periods. The maximum total costs over all $100$ demand vectors approximates the worst case costs with the FH approach, and is depicted in Figure \ref{fig:inv-worstcase}. From this picture it becomes clear that the optimal order quantity for the first time period is approximately $2.3$, which has a worst case performance of $10.8$.

We have solved the model for the full time horizon with the RC, with the AARC (where $c_1$ and $c_2$ are adjustable on the full $\vec{d}$, and $q_2$ is adjustable on $d_1$), and as a certain problem with $\vec{d} = \vec{5}$. The nominal solution gives $q_1=0$, the RC gives $q_1 \approx 4.4$, while the AARC yields $q_1 \approx 5.3$, leading to worst case costs of the FH approach of $17.8$, $14.9$ and $16.8$, respectively. So, the RC takes the best ``here and now'' decision with respect to the worst case performance. It may be paradoxical that the AARC yields a worse solution than the RC, since the feasible region of the AARC includes the RC solution. However, neither of the two optimize the right objective function. Both approximate the objective value using (static or adjustable) auxiliary variables $c_t$. While AARC indeed has a better objective value than RC, the solution is not better for the original objective function.

We also perform a comparison on a more realistic problem, which is the retailer-supplier flexible commitment problem by \citet{BenTal:2005:retailer}. At the starting time, the retailer commits himself to ordering certain quantities in later months. These commitments are flexible, i.e., deviations are allowed at a penalty cost. The objective is to minimize the total costs for the retailer, consisting of ordering costs (minus the salvage value), holding costs, backlogging costs, penalty costs associated with deviating from the commitments, and costs for variability in the commitments. The article provides two data sets for twelve time periods, A12 and W12, which we also use in our optimization and comparison.

In this problem the retailer faces an uncertain demand. Following \citet{BenTal:2005:retailer} we consider box uncertainty where the demand may deviate up to $\rho\%$ around the nominal value. For the simulation we draw demand vectors uniformly from this box region. For these demand vectors the nominal solution, RC and AARC are carried out completely. For the FH approach, the reoptimization is performed after each time period based on previously observed demand. 500 demand vectors were generated for each data set and each uncertainty level $\rho$, and the same demand vectors were used for all models. In addition, the PH solution was computed for each of these demand vectors.

The simulation results are listed in Tables \ref{tbl:simulatedretailer} and \ref{tbl:simulatedretailerPoU}. The PH solution always performs best, but cannot be used in practice, and hence, we ignore it for now. For data set A12, the nominal solution with FH results in the lowest average costs. This means that the nominal solution takes better ``here and now'' decisions than RC and AARC. Moreover, the RC-FH has lower average costs than the AARC-FH, so also the RC takes better ``here and now'' decisions than AARC. The advantage of the nominal FH solution compared to RC-FH and AARC-FH increases when the uncertainty set becomes larger. The same conclusions can be drawn for data set W12, except that the nominal solution is the best solution and FH leads to higher mean costs, and except that AARC (without folding horizon) outperforms RC. The Nominal-FH in Table \ref{tbl:simulatedretailer} is an example of the statement in Section \ref{B2} that a high standard deviation is not necessarily bad. Its standard deviation is higher than that of AARC-FH but this is explained by the high standard deviation of PH. In Table \ref{tbl:simulatedretailerPoU} we see that Nominal-FH has a lower mean and standard deviation than AARC-FH, which means that Nominal-FH is closer to the PH solution than AARC-FH. For data set A12 and $\rho=10\%$, it is not clear whether RC outperforms AARC. We now demonstrate the two statistical tests from Section \ref{B2} on this data set, each based on 100 newly generated uncertainty vectors, to test whether RC outperforms AARC. The null hypothesis that both solutions perform equally on average is rejected ($p=4.7 \cdot 10^{-12}$), and also the null hypothesis that the medians of the RC and AARC are equal is rejected ($p=2.8 \cdot 10^{-8}$). These results show that the AARC is not necessarily better than the RC and support the statement in Section \ref{B2} that a simulation is required for comparing solutions. As mentioned in Section \ref{sec:solvingrc}, RO may provide solutions that are not Pareto efficient when multiple optimal solutions exist. A different optimal solution to the RC or AARC may yield completely different simulation results, rendering our conclusions useless. We have therefore generated solutions that are robust optimal, and additionally, are optimal with respect to the nominal demand trajectory. Although these solutions were different, the conclusions remain unaffected.

\begin{table}
	\centering
  \caption{Simulated mean (std) costs for the retailer-supplier flexible commitment problem.}
  \label{tbl:simulatedretailer}
	\begin{tabular}{l ll ll}
	\toprule
	         &   \multicolumn{2}{c}{A12} & \multicolumn{2}{c}{W12}  \\
                \cmidrule(lr){2-3}    \cmidrule(lr){4-5}
           &   $\rho=10\%$ & $\rho=50\%$ &   $\rho=10\%$ & $\rho=50\%$ \\
   \midrule
Nominal    &  816 (35) &  999 (217) & 12799 (742) & 16002 (3634) \\
RC         &  861 (16) & 1080  (81) & 13907 (275) & 21564 (1415) \\
AARC       &  874 (4)  & 1030 (105) & 13313  (35) & 18577  (179) \\
\midrule
Nominal-FH &  807 (21) &  967 (156) & 12839 (288) & 16276 (1447) \\
RC-FH      &  829 (6)  & 1030  (47) & 13110  (56) & 17578 (309) \\
AARC-FH    &  881 (4)  & 1063  (65) & 13312  (35) & 18577 (178) \\
\midrule
PH         &  785 (18) &  826  (95) & 12187 (204) & 12950 (1024) \\
  \bottomrule
	\end{tabular}
\end{table}

\begin{table}
	\centering
  \caption{Simulated mean (std) PoU for the retailer-supplier flexible commitment problem.}
  \label{tbl:simulatedretailerPoU}
	\begin{tabular}{l ll ll}
	\toprule
	         &   \multicolumn{2}{c}{A12} & \multicolumn{2}{c}{W12}  \\
                \cmidrule(lr){2-3}    \cmidrule(lr){4-5}
           &   $\rho=10\%$ & $\rho=50\%$ &   $\rho=10\%$ & $\rho=50\%$ \\
   \midrule
Nominal    &   30 (25) & 173 (156) &  653 (597) & 3053 (2979) \\
RC         &   76 (29) & 253 (141) & 1720 (466) & 8615 (2370) \\
AARC       &   89 (16) & 228  (78) & 1126 (231) & 5627 (1160) \\
\midrule
Nominal-FH &   22  (9) & 140  (75) &  653 (144) & 3326  (703) \\
RC-FH      &   44 (16) & 204 (103) &  924 (207) & 4628 (1042) \\
AARC-FH    &   95 (16) & 237 (125) & 1125 (232) & 5527 (1161) \\
  \bottomrule
	\end{tabular}
\end{table}
\newpage

\section{Conclusion}\label{conc}
In this paper, we give a concise introduction and a general procedure that shall be helpful for using RO in practice. Additionally, we give several practical insights and hints in applying RO. Examples of such practical insights are: the robust reformulations of equivalent deterministic optimization problems may not be equivalent; in multi-stage optimization problems, re-optimizing the given problem at each stage using static RO or nominal data may outperform solutions provided by ARO; and the actual probability guarantee of an uncertainty set is often higher than the probabilistic guarantee that is approximated by using a safe approximation technique. We also discuss many practical issues to apply RO in a successful and convincing way. Examples are: how to choose the uncertainty set; what is the right interpretation of ``RO optimizes for the worst case''; and should the decision rule used in ARO be a function of the final or the primitive uncertainty? Moreover, we propose ideas on how to deal with equality constraints and integer adjustable variables, and on how to compare the robustness characteristics of two solutions. We have provided many numerical examples to illustrate our insights and discussions.

\bibliographystyle{abbrvnatnew}
\bibliography{libraryprac}

\end{document}